\documentclass{amsart}

\usepackage{amsmath,amssymb,amscd,amsthm,latexsym}
\usepackage[applemac]{inputenc}
\usepackage{graphics}
\usepackage[all]{xy}

\newtheorem{theor}{Theorem}[section]
\newtheorem{lem}[theor]{Lemma}
\newtheorem{prop}[theor]{Proposition}
\newtheorem{cor}[theor]{Corollary}
\newtheorem{defin}[theor]{Definition}
\newtheorem{rem}[theor]{Remark}

\newcommand{\dem}{\noindent \textit{Proof. }}
\newcommand{\findem}{\hfill $\Box$\\}

\newdir{ >}{{}*!/-7pt/\dir{>}}

\def\const{\mbox{\sf c}}

\title{Relative directed homotopy theory of partially ordered spaces}

\author{Thomas Kahl}
\begin{document}

\begin{abstract}
Algebraic topological methods have been used successfully in concurrency theory, the domain of theoretical computer science that deals with distributed computing. L. Fajstrup, E. Goubault, and M. Raussen have introduced partially ordered spaces (pospaces) as a model for concurrent systems. In this paper it is shown that the category of pospaces under a fixed pospace is both a fibration and a cofibration category in the sense of H. Baues. The homotopy notion in this fibration and cofibration category is relative directed homotopy. It is also shown that the category of pospaces is a closed model category such that the homotopy notion is directed homotopy. 
\end{abstract}

\maketitle


\section{Introduction}
It has turned out in the recent past that homotopy theoretical methods can be employed efficiently to study problems in concurrency theory. This is the domain of theoretical computer science that deals with parallel computing and distributed databases. Various topological models have been introduced in order to describe concurrent systems. Examples are partially ordered spaces \cite{FajstrupGR}, flows \cite{Gaucher}, globular CW-complexes \cite{GaucherGoubault}, and d-spaces \cite{Grandis}. The reader is referred to E. Goubault \cite{Goubault} for a recent introduction to different topological models for concurrency.

In this paper we shall study the homotopy theory of partially ordered spaces which have been introduced by  
L. Fajstrup, E. Goubault, and M. Raussen in \cite{FajstrupGR}. A partially ordered space (or pospace) is a topological space $X$ equipped with a partial order $\leq$. The space $X$ is interpreted as the state space of a concurrent system. The partial order $\leq$ represents the time flow. The idea here is that the execution of a system is a process in time so that a system in each state $x$ can only proceed to subsequent states $y \geq x$ and not go back to preceding states $y < x$. A natural question is whether a system in a given state $x$ can reach another state $y$ or, in other words, whether there is an “execution path" from $x$ to $y$. Such problems can be formalized appropriately using the following notion of maps between pospaces. A dimap (short for directed map) from a pospace $(X,\leq)$ to a pospace $(Y,\leq)$ is a continuous map $f : X \to Y$  such that $x\leq y$ implies $f(x) \leq f(y)$. An execution path from a state $x$ of a pospace $(X,\leq)$ to a state $y$ can now formally be defined to be a dimap $f$ from the unit interval $I = [0,1]$ with the natural order to $(X,\leq)$ such that $f(0) = x$ and $f(1) = y$. 

Consider a very simple concurrent system where two processes $A$ and $B$ modify asynchronously a shared resource. This situation can be modeled by the pospace $(X,\leq\nolinebreak)$ where $X = (I\times I) \setminus (]\frac{1}{3}, \frac{2}{3}[ \times ]\frac{1}{3}, \frac{2}{3}[)$ and $\leq$ is the componentwise natural order. If in a state $(x,y) \in X$, $x < \frac{1}{3}$ then $A$ has not yet accessed the resource; if $x = \frac{1}{3}$, $A$ has accessed the resource and is ready to modify it, if $\frac{1}{3} < x < \frac{2}{3}$ then $A$ is modifying the resource, and if $x \geq \frac{2}{3}$ then $A$ has modified the resource. Similarly, $B$ has not yet accessed, has accessed, modifies, and has modified the resource if $y \in [0,\frac{1}{3}[$, $y = \frac{1}{3 }$, $y\in ]\frac{1}{3},\frac{2}{3}[$, and $y \in [\frac{2}{3},1]$ respectively. Since the processes cannot modify the resource simultaneously, there are no possible states in $]\frac{1}{3}, \frac{2}{3}[\times ]\frac{1}{3}, \frac{2}{3}[$. The system has an initial state $(0,0)$ and a final state $(1,1)$ and there are infinitely many execution paths from $(0,0)$ to $(1,1)$. There are two kinds of such paths: those whose second coordinate is in $[0,\frac{1}{3}]$ when the first coordinate is in $]\frac{1}{3},\frac{2}{3}[$ and those whose second coordinate is in $[\frac{2}{3},1]$ when the first coordinate is in $]\frac{1}{3},\frac{2}{3}[$. The execution paths of the first kind correspond to executions where $A$ modifies the resource before $B$ and the execution paths of the second kind correspond to executions where $B$ modifies the resource before $A$. From a computer scientific point of view it makes therefore sense to regard execution paths of the same kind as equivalent. The equivalence relation behind this is dihomotopy (short for directed homotopy) relative to the initial and final states.
As the name suggests, this is a kind of homotopy and so homotopy theory becomes relevant to concurrency theory.

Before we define dihomotopy we note that for every topological space $X$ the diagonal $\Delta \subset X\times X$ is a partial order. We also note that the product of two pospaces exists in the category theoretical sense and is the topological product with the component\-wise order. 
Two dimaps $f, g: (X\leq) \to (Y,\leq)$ are said to be dihomotopic if there exists a dimap $H : (X,\leq) \times (I,\Delta) \to (Y,\leq)$ such that $H(x,0) = f(x)$ and $H(x,1) = g(x)$. The example above shows that one also needs a relative notion of dihomotopy. Indeed, in the absolute sense, any execution path is dihomotopic to a constant dimap.
As P. Bubenik \cite{Bubenik} has pointed out, another reason for considering a relative notion of dihomotopy is the fact that it depends a lot on the context whether two pospaces can be interpreted as models of the same concurrent system. In order to define relative dihomotopy we work in the comma category of pospaces under a fixed pospace $(C,\leq)$. A $(C,\leq)$-pospace is a triple $(X, \leq, \xi)$ consisting of a pospace $(X,\leq)$ and a dimap $\xi : (C,\leq) \to (X,\leq)$. A $(C,\leq)$-dimap $f:(X,\leq,\xi) \to (Y,\leq,\theta)$ is a dimap $f: (X,\leq) \to (Y,\leq)$ such that $f \circ \xi = \theta$. Two $(C,\leq)$-dimaps $f,g : (X,\leq,\xi) \to (Y,\leq,\theta)$ are said to be dihomotopic relative to $(C,\leq)$ if there exists a dimap $H: (X,\leq)\times (I,\Delta) \to (Y,\leq)$ such that $H(x,0) = f(x)$, $H(x,1) = g(x)$ $(x\in X)$, and $H(\xi(c),t) = \theta (c)$ $(c \in C$, $t\in I)$. In the above example let $(C,\leq)$ be the discrete space $\{0,1\}$ with the natural order and consider the inclusion $\iota : \{0,1\} \hookrightarrow I$ and the dimap $\xi :(\{0,1\},\leq) \to (X,\leq)$ given by  $\xi(0) = (0,0)$ and $\xi(1) = (1,1)$. Then two execution paths from $(0,0)$ to $(1,1)$, i.e., two $(\{0,1\},\leq )$-dimaps $f,g : (I,\leq,\iota) \to (X,\leq,\xi)$, are of the same kind if and only if they are dihomotopic relative to $(\{0,1\},\leq)$.

The best known framework for homotopy theory is certainly the one of closed model categories in the sense of D. Quillen \cite{Quillen}. A closed model category is a category with three classes of morphisms, called weak equivalences, fibrations, and cofibrations, which are subject to certain axioms. The structure of a closed model category splits up into two dual structures which are essentially the structure of a cofibration category and the structure of a fibration category. Cofibration and fibration categories have been introduced by H. Baues \cite{Baues} who has developed an extensive homotopy theory for these categories. In this paper we show that the category of $(C,\leq)$-pospaces is both a fibration and a cofibration category (Theorems \ref{fibmain} and \ref{pocofcat}). We also show that the category of absolute pospaces is a closed model category (Theorem \ref{model}). The main ingredient of the homotopy theory of a cofibration, fibration, or closed model category is of course a notion of homotopy. We show that this notion of homotopy in the cofibration and fibration category of $(C,\leq)$-pospaces is dihomotopy relative to $(C,\leq)$ (cf. \ref{hoverB} and \ref{pocofcat}). Similarly, the homotopy notion of the closed model category of pospaces is dihomotopy (cf. \ref{model}). The reader should note that the model structure on the category of pospaces induces a model structure on the category of $(C,\leq)$-pospaces. This means that some relative dihomotopy theory is part of absolute dihomotopy theory. In contrast to what is happening in ordinary homotopy theory, the relative dihomotopy theory coming from the model structure on absolute pospaces is too restrictive. The reason for this is that one cannot restrict oneself to $(C,\leq)$-pospaces $(X,\leq,\xi)$ where the dimap $\xi$ is a dicofibration, i.e., a cofibration in the closed model category of pospaces. It is therefore necessary to develop relative dihomotopy theory in an autonomous way.

Pospaces are a rather simple model for concurrency. L. Fajstrup, E. Goubault, and M. Raussen \cite{FajstrupGR} also introduce locally partially ordered spaces, or local pospaces, which consitute a more advanced model for concurrency. There are dimaps of local pospaces and there is a concept of dihomotopy. One can show that the category of local pospaces is a fibration category such that the homotopy notion is dihomotopy. It is, however, not known whether there are enough colimits for a cofibration or a closed model category structure. Note that P. Bubenik and K. Worytkiewicz \cite{BubenikW} have constructed a closed model category containing the category of local pospaces as a subcategory. Another more sophisticated model for concurrency is given by flows \cite{Gaucher}. Unfortunately, there seem to be serious problems in constructing a model structure for dihomotopy on the category of flows (cf. \cite{GaucherNoModel}). Pospaces are probably not suited to all aspects of concurrency but they allow a rather straightforward extension of ordinary homotopy theory. 

The paper is organized as follows. In section \ref{Pospaces} we show that the category of $(C,\leq\nolinebreak)$-pospaces is complete and cocomplete. Section \ref{dihomotopy} contains the fundamental material about dihomotopy. In particular, we define dihomotopy equivalences relative to $(C,\leq)$ and the adjoint cylinder and path $(C,\leq)$-pospace functors. In section \ref{difibrations} we define $(C,\leq)$-difibrations and prove some fundamental facts about them. 
The main result of section \ref{fibrationcategory} is Theorem \ref{fibmain} which states that the category of $(C,\leq)$-pospaces is a 
fibration category where fibrations are $(C,\leq)$-difibrations and weak equivalences are dihomotopy equivalences relative to $(C,\leq)$. This result is a consequence of the fact that the the category of $(C,\leq)$-pospaces is a P-category in the sense of \cite{Baues} which is proved in \ref{PpoTop}. Proposition \ref{hoverB} contains the result that two $(C,\leq)$-dimaps are homotopic in the fibration category of $(C,\leq)$-pospaces if and only if they are dihomotopic relative to $(C,\leq)$. In section \ref{cofibrations} we study cofibrations in a fibration category and show in Theorem \ref{cofibrationcategory} that they induce under certain conditions the structure of a cofibration category. We show that the homotopy notions of the cofibration and fibration category structures coincide (cf. \ref{ho=ho}). The internal cofibrations of the fibration category of $(C,\leq)$-pospaces are called $(C,\leq)$-dicofibrations. In \ref{pocofcat} we show that the conditions of \ref{cofibrationcategory} are satisfied so that the category of $(C,\leq)$-pospaces is a cofibration category in which the homotopy notion is dihomotopy relative to $(C,\leq)$. In the last section it is shown that the category of absolute pospaces is a closed model category such that the homotopy notion is dihomotopy.

\section{Pospaces} \label{Pospaces}

\begin{defin}\rm
A \textit{pospace} (short for \textit{partially ordered space}) is a pair $(X, \leq)$ consisting of a
space $X$ and a partial order $\leq$ on $X$. If this is helpful we shall write $\leq_X$ instead of $\leq$. A \textit{dimap} (short for \textit{directed map}) $f : (X, \leq) \to (Y, \leq)$ is a continuous map $f : X \to Y$ such that for all $x,x' \in X$, $x \leq x'$ implies $f(x) \leq f(x')$. The category of pospaces will be denoted by $\bf poTop$.
\end{defin}

In the original definition (cf. \cite{FajstrupGR}) the partial order of a pospace $(X,\leq)$ is required to be closed as a subspace of $X \times X$. It is shown in \cite{FajstrupGR} that a space $X$ can be equipped with such a closed partial order if and only if it is a Hausdorff space, and in some sense pospaces with a closed partial order are for general pospaces what Hausdorff spaces are for general topological spaces. An interesting topological space will of course in general be a Hausdorff space. From the homotopy theoretical point of view, however, a restriction to Hausdorff spaces is not necessary and it is indeed easier to develop ordinary homotopy theory in the category of all topological spaces than in the category of Hausdorff spaces. For the same reason of simplicity we shall work with general pospaces rather than with pospaces having a closed partial order.   

For every topological space $X$ the diagonal $\Delta \subset X \times X$ is a partial order. The functor $X \mapsto (X,\Delta)$ from the category $\bf Top$ of topological spaces to $\bf poTop$ is left adjoint to the forgetful functor ${\bf poTop} \to {\bf Top}$.

\begin{prop} \label{complete}
The category $\bf poTop$ is complete and cocomplete. 
\end{prop}

\dem
We show that $\bf poTop$ has arbitrary products, coproducts, equalizers, and coequalizers. 
Let $\{ (X_i,\leq_{X_i})\}_{i\in I}$ be a family of pospaces.
The product of the pospaces $(X_i,\leq_{X_i})$ is the topological product $\prod _{i\in I} X_i$ equipped with the componentwise partial order. The coproduct of the pospaces $(X_i,\leq_{X_i})$ is the topological coproduct $\coprod _{i\in I} X_i$ together with the partial order given by 
$$x \leq y \Leftrightarrow \exists \, i\in I : x, y \in X_i \; \mbox{and}\; x\leq_{X_i} y.$$
Let $f, g: (X,\leq) \to (Y,\leq)$ be two dimaps. The equalizer of $f$ and $g$ is the topological equalizer $\{ x \in X\,|\, f(x) = g(x)\}$ together with the restriction of $\leq_X$. 

Consider the coequalizer of $f$ and $g$ in in $\bf Top$. This is the quotient space $Y/\sim$ where $\sim$ is the equivalence relation given by $f(x) \sim g(x)$ (i.e., the equivalence relation induced by the relation $f(x) \sim g(x)$). Define a relation $\vartriangleleft$ on $Y/\sim$ by 
$$\alpha \vartriangleleft \beta \Leftrightarrow \exists \, y_1, \dots , y_n \in Y : y_1 \in \alpha, y_n \in \beta,\;  \mbox{and}\; y_1 \leq y_2 \sim y_3 \leq \cdots \sim y_{n-1} \leq y_n.$$
The relation $\vartriangleleft$ is reflexive and transitive but not necessarily antisymmetric. Consider the equivalence relation $\vartriangleleft \vartriangleright$ on $Y/\sim$ defined by 
$$\alpha \vartriangleleft \vartriangleright \beta \Leftrightarrow \alpha \vartriangleleft \beta \; \mbox{and} \; \beta \vartriangleleft \alpha.$$
Define a relation $\leq$ on the quotient space $(Y/\sim)/\vartriangleleft \vartriangleright$ by 
$$A \leq B \Leftrightarrow \forall \, \alpha \in A, \beta \in B : \alpha \vartriangleleft \beta.$$
This is a partial order. Let $p : Y \to (Y/\sim)/\vartriangleleft \vartriangleright$ be the identification map and consider two elements $y \leq z \in Y$. Then $[y] \vartriangleleft [z]$. Let $\alpha \in [[y]] = p(y)$ and $\beta \in [[z]] = p(z)$. Then $\alpha \vartriangleleft [y] \vartriangleleft [z] \vartriangleleft \beta$ and hence $\alpha \vartriangleleft \beta $. This implies that $p$ is a dimap. We show that the pospace $((Y/\sim)/\vartriangleleft \vartriangleright, \leq)$ together with the dimap $p$ is the coequalizer of $f$ and $g$ in $\bf poTop$. Let $h : (Y,\leq) \to (Z,\leq)$ be a dimap such that $h \circ f = h\circ g$. Let $\bar h : Y/\sim \to Z$ be the continuous map indued by $f$ and $g$, i.e., $\bar h ([y]) = h(y)$. Let $\alpha , \beta \in Y/\sim $ such that $\alpha \vartriangleleft \beta$. Then there exist $y_1, \dots , y_n \in Y$ such that $y_1 \in \alpha$, $y_n \in \beta$, and  $$y_1 \leq y_2 \sim y_3 \leq \cdots \sim y_{n-1} \leq y_n.$$
It follows that  
$$\bar h (\alpha) = h(y_1) \leq h(y_2) = h(y_3) \leq \cdots = h(y_{n-1}) \leq h(y_n) = \bar h(\beta)$$
and hence that $\bar h(\alpha) \leq \bar h(\beta)$. If $\alpha \vartriangleleft \vartriangleright \beta$ we therefore have $\bar h(\alpha) = \bar h (\beta)$. There hence exists a unique continuous map $\tilde h : (Y/\sim)/\vartriangleleft \vartriangleright \to Z$ such that $\tilde h ([\alpha]) = \bar h(\alpha)$. We have $\tilde h (p(y)) = \bar h([y]) = h(y)$. Let $\alpha, \beta \in Y/\sim$ such that $[\alpha] \leq [\beta]$. Then 
$\alpha \vartriangleleft \beta $ and  $\tilde h([\alpha]) = \bar h(\alpha ) \leq \bar h(\beta) = \tilde h([\beta])$.
It follows that $\tilde h$ is a dimap $((Y/\sim)/\vartriangleleft \vartriangleright,\leq \nolinebreak) \linebreak \to (Z,\leq)$ satisfying $\tilde h \circ p = h$. Since $p$ is surjective, $\tilde h$ is the only dimap with this property.
\findem

\begin{defin}\rm
Let $(C,\leq)$ be a pospace. A \textit{$(C,\leq)$-pospace} is a triple $(X,\leq,\xi)$ consisting of a pospace $(X,\leq)$ and a dimap $\xi : (C,\leq) \to (X,\leq)$. A \textit{$(C,\leq)$-dimap} from $(X,\leq,\xi)$ to $(Y,\leq,\theta)$ is a dimap $f:(X,\leq) \to (Y,\leq)$ such that $f\circ \xi = \theta$. The category of $(C,\leq)$-pospaces is denoted by $(C,\leq)$-${\bf poTop}$.
\end{defin}

\begin{prop} \label{ccomplete}
For any pospace $(C,\leq)$ the category $(C,\leq)$-${\bf poTop}$ is complete and cocomplete. 
\end{prop}

\dem
This follows from \ref{complete}.
\findem

\begin{rem}\rm
An absolute pospace is the same as a $(\emptyset, \Delta)$-pospace.
\end{rem}

\section{Relative Dihomotopy} \label{dihomotopy}
Throughout this section we work under a fixed pospace $(C,\leq)$. We define dihomotopy relative to $(C,\leq)$, introduce the adjoint cylinder and path $(C,\leq)$-pospace functors, and give characterizations of relative dihomotopy by means of these constructions.

\begin{defin}\rm
Two $(C,\leq)$-dimaps $f,g : (X,\leq,\xi) \to (Y, \leq,\theta)$ are said to be \textit{dihomotopic relative to $(C,\leq)$},
$f \simeq g \; rel. \; (C,\leq)$, if there exists a \textit{dihomotopy relative to $(C,\leq)$} from $f$ to $g$, i.e., a dimap $H : (X,\leq) \times (I,\Delta ) \to (Y,\leq)$ such that $H(x,0) = f(x)$, $H(x,1) = g(x)$ $(x\in X)$, and $H(\xi(c),t) = \theta(c)$ $(c\in C$, $t \in I)$. If $C = \emptyset$ we simply talk of dihomotopies and dihomotopic dimaps and we simply write $f\simeq g$.
\end{defin}

\begin{prop}
Dihomotopy relative to $(C,\leq)$ is an equivalence relation on the set of $(C,\leq)$-dimaps from
$(X,\leq,\xi)$ to $(Y, \leq,\theta)$. Furthermore, dihomotopy relative to $(C,\leq)$ is compatible with composition.
\end{prop}

\dem
This is an easy exercise.
\findem

\begin{defin}\rm
The equivalence class of a $(C,\leq)$-dimap with respect to dihomotopy relative to $(C,\leq)$ is called its \textit{dihomotopy class relative to $(C,\leq)$}. The quotient category $(C,\leq)$-${\bf poTop}/\simeq \; rel. \; (C,\leq)$ is the \textit{dihomotopy category relative to $(C,\leq)$}. A \textit{dihomotopy equivalence relative to $(C,\leq)$} is a $(C,\leq)$-dimap \linebreak $f : (X, \leq,\xi) \to (Y, \leq,\theta)$ such
that there exists a \textit{dihomotopy inverse relative to $(C,\leq)$}, i.e., a $(C,\leq)$-dimap $g : (Y,\leq,\theta) \to (X, \leq,\xi)$
satisfying $f\circ g \simeq id_{(Y,\leq,\theta)}$ $rel. \; (C,\leq)$ and $g \circ f \simeq id_{(X,\leq,\xi)}\; rel.\; (C,\leq)$. Two $(C,\leq)$-pospaces $(X,\leq,\xi)$ and
$(Y,\leq,\theta)$ are said to be \textit{dihomotopy equivalent relative to $(C,\leq)$} or of the \textit{same dihomotopy type relative to $(C,\leq)$} if there exists a dihomotopy equivalence relative to $(C,\leq)$ from $(X, \leq,\xi)$ to $(Y, \leq,\theta)$. If $C = \emptyset$ we simply talk of dihomotopy classes, the dihomotopy category, dihomotopy equivalences, and dihomotopy equivalent pospaces.
\end{defin}

Note that a $(C,\leq)$-dimap  is a dihomotopy equivalence relative to $(C,\leq)$ if and only if its dihomotopy class relative to $(C,\leq)$ is an
isomorphism in the dihomotopy category relative to $(C,\leq)$. Similarly, two $(C,\leq)$-pospaces are dihomotopy equivalent relative to $(C,\leq)$  
if and only if they are isomorphic in the dihomotopy category relative to $(C,\leq)$.

\begin{prop} \label{2=3}
Any isomorphism of $(C,\leq)$-pospaces is a dihomotopy equivalence relative to $(C,\leq)$. Let $f: (X, \leq,\xi) \to (Y,\leq,\theta)$ and $g:(Y, \leq,\theta) \to (Z, \leq,\zeta )$ be two $(C,\leq)$-dimaps. If two of
$f$, $g$, and $g\circ f$ are dihomotopy equivalences relative to $(C,\leq)$, so is the third. Any retract of a dihomotopy equivalence relative to $(C,\leq)$ is a dihomotopy equivalence relative to $(C,\leq)$. 
\end{prop}

\dem
The first statement is obvious and the others follow from the corresponding facts for isomorphisms.
\findem

Let $(X,\leq,\xi)$ be a $(C,\leq)$-pospace and $S$ be a space. Form the pushout diagram of pospaces
$$\xymatrix{
(C,\leq) \times (S,\Delta ) \ar[r]^-{pr_C}_{} \ar[d]_{\xi \times id_S}_{} 
& (C,\leq) \ar [d]^{\bar \xi} 
\\ 
(X,\leq) \times (S,\Delta) \ar[r]^{}_{} 
& (X\Box_C S,\leq).
}$$ 
The space $X \Box_C S$ is the pushout of the underlying diagram of spaces. If $S = \emptyset$, $(X\Box_C S,\leq) = (C,\leq)$. If $S\not= \emptyset$, we may construct $X\Box_C S$ as the quotient space $(X\times S)/\sim$ where $$(x,s)\sim (y,t) \Leftrightarrow (x,s) = (y,t)   \; \mbox{or}\; \exists \, c \in C : x = y = \xi(c).$$ The partial order on $X\Box_C S$ is then given by 
$$[x,s] \leq [x',s'] \Leftrightarrow (x,s)\leq (x',s')  \; \mbox{or}\;  \exists \, c \in C : x \leq \xi (c) \leq x'.$$ 
We define a pospace under $(C,\leq)$ by setting $$(X,\leq,\xi) \Box_{(C,\leq)} S = (X\Box_C S, \leq, \bar \xi).$$ It is clear that this construction is natural and
defines a functor \begin{center}{$\Box_{(C,\leq)} : (C,\leq)$-${\bf poTop} \times {\bf Top} \to (C,\leq)$-${\bf poTop}$.}\end{center}

\begin{defin}\rm \label{cylinder}
The \textit{cylinder} on a $(C,\leq)$-pospace $(X,\leq,\xi)$ is the $(C,\leq)$-pospace $(X,\leq,\xi) \Box_{(C,\leq)} I$. 
\end{defin}

Note that if $C = \emptyset$ then the cylinder on a pospace $(X,\leq)$ is just the product pospace $(X,\leq) \times (I,\Delta)$.

\begin{prop}
Two $(C,\leq)$-dimaps $f, g: (X,\leq, \xi) \to (Y,\leq,\theta)$ are dihomotopic relative to $(C,\leq)$ if and only if there exists a $(C,\leq)$-dimap $H: (X,\leq,\xi)\Box_{(C,\leq)} I \to \linebreak (Y,\leq,\theta)$ such that $H([x,0]) = f(x)$ and $H([x,1]) = g(x)$.
\end{prop}

\dem
This is straightforward.
\findem

Recall that the path space $X^I$ of a topological space $X$ is the set of all continuous maps $\omega : I \to X$ with the compact-open topology.

\begin{defin}\rm \label{pathpospace}
Let $(X ,\leq , \xi)$ be a $(C,\leq)$-pospace. The \textit{path $(C,\leq)$-pospace} of $(X,\leq,\xi)$ is the $(C,\leq)$-pospace $(X^I,\leq,\const_{\xi})$ where the partial order is given by
$$\omega \leq \nu \Leftrightarrow \forall \, t \in I:  \omega (t) \leq \nu (t)$$
and the dimap $\const_{\xi} : (C,\leq) \to (X^I, \leq)$ is given by  $\const_{\xi (c)}(t) = \xi(c)$ $(c\in C$, $t \in I)$. 
\end{defin}

The path $(C,\leq)$-pospace is obviously functorial. Note also that for each $t\in I$ the evaluation map $ev_t : X^I \to X$, $\omega \mapsto \omega (t)$ is a $(C,\leq)$-dimap $(X^I,\leq,\const_{\xi}) \to (X ,\leq , \xi)$.

\begin{prop} \label{adjoints}
The path $(C,\leq)$-pospace functor is right adjoint to the cylinder functor $-\Box_{(C,\leq)}I$.
\end{prop}

\dem
The natural correspondence between the $(C,\leq)$-dimaps $h : (X,\leq,\xi) \to (Y^I,\leq, \const_{\theta})$ and the $(C,\leq)$-dimaps $H : (X,\leq,\xi) \Box_{(C,\leq)}I \to (Y,\leq,\theta)$ is given by the formula $h(x)(t) = H([x,t])$.
\findem

Using this adjunction one easily establishes the following characterization of dihomotopy relative to $(C,\leq)$:

\begin{prop}  \label{pathhomotopic}
Two $(C,\leq)$-dimaps $f,g : (X, \leq,\xi) \to (Y, \leq,\theta)$ are dihomotopic relative to $(C,\leq)$ if and only if there exists a
$(C,\leq)$-dimap $h : (X, \leq,\xi) \to \linebreak (Y^I , \leq,\const_{\theta})$ such that $f = ev_0\circ h$ and $g = ev_1 \circ h$.
\end{prop}

\section{$(C,\leq)$-Difibrations} \label{difibrations}

As in the preceding section we work under a fixed pospace $(C,\leq)$. Recall the following terminology:

\begin{defin}\rm 
Let $\bf C$ be a category and $\mathcal A$ be a class of morphisms. A morphism $f: X \to Y$ is said to have the \textit{right lifting property} with respect to $\mathcal A$ if for every morphism $a : A \to B$ of $\mathcal A$ and for all morphisms $g : A \to X$ and $h: B \to Y$ satisfying $f \circ g = h\circ a$ there exists a morphism $\lambda : B \to X$ such that $f \circ \lambda = h$ and $\lambda \circ a = g$. Similarly, a morphism $f: X \to Y$ is said to have the \textit{left lifting property} with respect to $\mathcal A$ if for every morphism $a : A \to B$ of $\mathcal A$ and for all morphisms $g : X \to A$ and $h: Y \to B$ satisfying $a \circ g = h\circ f$ there exists a morphism $\lambda : Y \to A$ such that $a \circ \lambda = h$ and $\lambda \circ f = g$. 
\end{defin}

\begin{defin} \rm
A \textit{$(C,\leq)$-difibration} is a $(C,\leq)$-dimap having the right lifting property with respect to the  $(C,\leq)$-dimaps of the
form $$i_0 : (X,\leq,\xi) \to (X,\leq,\xi)\Box_{(C,\leq)}  I,\quad i_0(x) = [x,0].$$
If $C = \emptyset$ we simply talk of difibrations.
\end{defin}

It is a general fact that any class of morphisms in a category which is defined by having the right (resp. left) lifting property with respect to another class of morphisms contains all isomorphisms and is closed under base change (resp. cobase change), composition, and retracts. We therefore have

\begin{prop} \label{pullbackfibration}
The class of $(C,\leq)$-difibrations is closed under composition, retracts, and base change. Every isomorphism of $(C,\leq)$-pospaces is a
$(C,\leq)$-difibration.
\end{prop}

We leave it to the reader to check the following $\Box$-free characterization of $(C,\leq)$-difibrations:

\begin{prop}
A $(C,\leq)$-dimap $p : (E, \leq,\varepsilon ) \to (B, \leq,\beta )$ is a $(C,\leq)$-difibration if and only if for every $(C,\leq)$-dimap $f : (X, \leq,\xi)  \to (E, \leq,\varepsilon)$ and every dimap \linebreak 
$H : (X,\leq )\times (I,\Delta) \to (B, \leq)$ satisfying $H(x,0) = (p\circ f)(x)$ $(x\in X)$ and $H(\xi(c),t) = \beta(c)$ $(c \in C)$ there exists a dimap 
$G : (X,\leq) \times (I,\Delta) \to (E, \leq)$ such that $G(x, 0) = f(x)$ $(x\in X)$, $p \circ G = H$, and $G(\xi(c),t) = \varepsilon(c)$ $(c\in C$, $t \in I)$.
\end{prop}

\begin{prop} \label{allobjectsfibrant}
For every $(C,\leq)$-pospace $(X, \leq,\xi)$ the final $(C,\leq)$-dimap \linebreak $\ast : (X, \leq,\xi) \to (\ast , \Delta ,\ast)$ is a $(C,\leq)$-difibration.
\end{prop}

\dem
Let $f : (W,\leq,\psi) \to (X, \leq\nolinebreak, \xi)$ be a $(C,\leq)$-dimap and
$F: (W,\leq) \times (I,\Delta) \to (\ast,\Delta )$ be a (the only) dimap. Define a  dimap
$H :  (W,\leq)  \times (I,\Delta) \to (X,\leq )$ by $H(w,t) = f(w)$. Then $H(w,0) = f(w)$, $\ast \circ H = F$, and $H(\psi (c),t) =\xi(c)$.
\findem

It is a very useful fact in ordinary homotopy theory (due to A. Str\o m \cite{NC1}) that fibrations have a much stronger lifting property than the defining homotopy lifting property. The last point of this section is an adaptation of this result to $(C,\leq)$-difibrations. We shall need the following lemma:

\begin{lem}\label{Boxass}
Let $(X,\leq, \xi)$ be a pospace and $S$ be a space. Then 
$$(X,\leq,\xi) \Box_{(C,\leq)} (S\times I) = ((X,\leq,\xi) \Box_{(C,\leq)} S) \Box_{(C,\leq)} I.$$
\end{lem}

\dem
Consider the defining pushout of $(X,\leq,\xi) \Box_{(C,\leq)} S$:
$$\xymatrix{
(C,\leq) \times (S,\Delta ) \ar[r]^-{pr_C}_{} \ar[d]_{\xi \times id_S}_{} 
& (C,\leq) \ar [d]^{\bar \xi} 
\\ 
(X,\leq) \times (S,\Delta) \ar[r]^{}_{} 
& (X\Box_C S,\leq).
}$$ 
Since  the functor $-\times (I,\Delta) : {\bf poTop} \to {\bf poTop}$ is a left adjoint, it preserves colimits. It follows that both squares in the following diagram of pospaces are pushouts:
$$\xymatrix{
(C,\leq) \times (S,\Delta )\times (I,\Delta)  \ar[r]^-{pr_C\times id_I}_{} \ar[d]_{\xi \times id_S \times id_I}_{} 
& (C,\leq)\times (I,\Delta) \ar [d]^{\bar \xi \times id_I} \ar[r]^{pr_C}
& (C,\leq) \ar [d]^{\bar {\bar \xi}} 
\\ 
(X,\leq) \times (S,\Delta)\times (I,\Delta) \ar[r]^{}_{} 
& (X\Box_C S,\leq) \times (I,\Delta) \ar[r] 
& ((X\Box_C S)\Box_C I,\leq).
}$$ 
This implies that the whole diagram is the defining pushout of $(X,\leq,\xi) \Box_{(C,\leq)} (S\times I)$ and thus that $(X,\leq,\xi) \Box_{(C,\leq)} (S\times I) = ((X,\leq,\xi) \Box_{(C,\leq)} S) \Box_{(C,\leq)} I$.
\findem

Recall that a \textit{trivial cofibration} of spaces is a closed cofibration which also is a homotopy equivalence. The following characterization of $(C,\leq)$-difibrations is of fundamental importance:

\begin{prop}\label{RLPBox}
A $(C,\leq)$-dimap $p : (E,\leq,\varepsilon) \to (B, \leq,\beta)$ is a $(C,\leq)$-difibration if and only if for every
$(C,\leq)$-pospace $(Z,\leq,\zeta)$, every trivial cofibration of spaces $i : A \to X$, every dimap
$f : (Z,\leq) \times (A,\Delta) \to (E,\leq\nolinebreak)$ satisfying $f(\zeta (c),a) = \varepsilon (c)$ $(c\in C$, $a\in A)$, and every dimap
$g : (Z,\leq) \times (X,\Delta) \to (B,\leq)$ satisfying $g(z,i(a))= p(f(z,a))$ $(z \in Z$, $a\in A)$ and $g(\zeta (c),x) = \beta (c)$ $(c \in C$, $x\in X)$ there exists a dimap
$\lambda : (Z,\leq \nolinebreak) \times \nolinebreak (\nolinebreak X,\Delta \nolinebreak) \to (E,\leq)$ such that
$\lambda (z,i(a))=  f(z,a)$ $(z \in Z$, $a\in A)$, $p\circ \lambda = g$, and $\lambda (\zeta (c),x) = \varepsilon (c)$ $(c \in C$, $x\in X)$.
\end{prop}

\dem
If $p$ has this lifting property, it is a $(C,\leq)$-difibration: it suffices to consider the trivial cofibration $\{0\} \hookrightarrow I$. Suppose that $p$ is a $(C,\leq)$-difibration and consider a $(C,\leq)$-pospace $(Z,\leq,\zeta)$, a trivial cofibration of spaces $i : A \to X$, a dimap $f : (Z,\leq) \times (A,\Delta) \to (E,\leq\nolinebreak)$ satisfying $f(\zeta (c),a) = \varepsilon (c)$, and a dimap \linebreak $g : (Z,\leq) \times (X,\Delta) \to (B,\leq)$ satisfying $g(z,i(a))= p(f(z,a))$ and $g(\zeta (c),x) = \beta (c)$. Since $i$ is a trivial cofibration, $i$ is a closed inclusion and $A$ is a strong deformation retract of $X$. There hence exist a retraction $r : X \to A$ of $i$ and a homotopy $H : X \times I \to X$ such that $H(x,0) = r(x)$, $H(x,1) = x$ $(x\in X)$, and $H(a,t) = a$ $(a\in A,\, t\in I)$. There also exists a continuous map $\phi : X \to I$ such that $A = \phi ^{-1}(0)$. Consider the map $G : X \times I \to X$ defined by
$$G(x,t)= \left\{ \begin{array}{lll}
H(x, \frac{t}{\phi(x)}) & t < \phi(x), \\
x & t \geq \phi(x).
\end{array}\right. $$
As in \cite[7.15]{Whitehead} one shows that $G$ is continuous. We have $G(x,0) = (i\circ r)(x)$ for all $x \in X$. Consider the following commutative diagram of
$(C,\leq)$-pospaces where $\bar f$ and $\bar g$ are given by $\bar f([z,x]) = f(z,r(x))$ and $\bar g([z,x,t]) = g(z,G(x,t))$:
$$\xymatrix{
(Z,\leq,\zeta)\Box_{(C,\leq)} X \ar[rr]^-{\bar f}_{} \ar[d]^{}_{id_Z\Box_C i_0}
& &(E,\leq,\varepsilon)\ar [d]^{p}_{}
\\
(Z,\leq,\zeta)\Box_{(C,\leq)} (X\times I) \ar[rr]^{}_-{\bar g}
& &(B,\leq,\beta).
}$$
By \ref{Boxass}, we may identify the $(C,\leq)$-dimap $id_Z\Box_Ci_0$ with the $(C,\leq)$-dimap $$(Z,\leq,\zeta)\Box _{(C,\leq)}X \to ((Z,\leq,\zeta)\Box _{(C,\leq)}X)\Box _{(C,\leq)} I,\quad [z,x] \mapsto [[z,x],0].$$
Since $p$ is a $(C,\leq)$-difibration, there exists a $(C,\leq)$-dimap
$$F : (Z,\leq,\zeta)\Box_{(C,\leq)} (X \times I) \to (E,\leq,\varepsilon)$$ such that $F \circ (id_Z\Box_C i_0) = \bar f$ and $p\circ F = \bar g$.
Consider the dimap \linebreak $\lambda : (Z,\leq) \times (X,\Delta ) \to (E,\leq)$ defined by $\lambda (z,x) = F([z,x,\phi(x)])$. We have
$$(p\circ \lambda) (z,x) = p(F([z,x,\phi(x)])) = \bar g([z,x,\phi(x)]) = g(z,G(x,\phi(x))) = g(z,x),$$
$$\lambda (z,i(a)) = \lambda (z,a) = F([z,a,\phi(a)]) = F([z,a,0]) = \bar f ([z,a]) = f(z,r(a)) = f(z,a),$$
and $\lambda (\zeta (c), x) = F([\zeta (c),x,\phi(x)]) = \varepsilon (c)$. This shows that $p$ has the required lifting property.
\findem

\section{The fibration category structure} \label{fibrationcategory}

The first result of this section is the fact that the category of $(C,\leq)$-pospaces is a \linebreak P-category in the sense of the following definition:

\begin{defin} \cite[I.3]{Baues} \rm
A category $\bf C$ equipped with a class of morphisms, called \textit{fibrations}, and a \textit{path object functor} $P: {\bf C} \to {\bf C}$, $X \mapsto X^I$, $f \mapsto f^I$ is said to be a \textit{P-category}
if it has a final object $\ast$ and if the following axioms are satisfied:
\begin{itemize}
\item[P1] There are natural transformations $q_0,q_1 : P \to id_{\bf C}$, $c: id_{\bf C} \to P$ such that $q_0\circ c = q_1\circ c = id$.
\item[P2] The pullback of two morphisms one of which is a fibration exists. The functor $P$ carries such a pullback into a pullback and
preserves the final object. The fibrations are closed under base change.
\item[P3] The composite of two fibrations is a fibration. Every isomorphism is a fibration and every final morphism $X \to \ast$ is a fibration. Every fibration $p: E \to B$ has the \textit{homotopy lifting property}, i.e., given morphisms $h: X \to B^I$ and $f : X \to E$ such that
$p\circ f = q_{\tau} \circ h$ $(\tau = 0$ or $\tau=1)$, there exists a morphism $H : X \to E^I$ such that $q_{\tau}\circ H = f$ and
$p^I \circ H = h$.
\item[P4] For every fibration $p: E \to B$ the morphism $$(q_0,q_1,p^I) : E^I \to (E\times E)\times _{B\times B}B^I$$ is a fibration.
Here, the target object is the fibered product of the morphisms $p\times p$ and $(q_0,q_1) : B^I \to B\times B$.
\item[P5] For each object $X$ there exists a morphism $T : (X^I)^I \to (X^I)^I$ such that $q_{\tau}^I\circ T = q_{\tau}$ and
$q_{\tau}\circ T = q_{\tau}^I$ $(\tau = 0,1)$.
\end{itemize}
\end{defin}

\begin{theor} \label{PpoTop}
Let $(C,\leq)$ be a pospace. The category $(C,\leq)$-${\bf poTop}$ is a P-cat\-egory. The fibrations are the $(C,\leq)$-difibrations and the functor $P$ is the path $(C,\leq)$-pospace functor.
\end{theor}

\dem
The natural transformations $q_0$ and $q_1$ are the evaluation maps $ev_0$ and $ev_1$. The natural transformation
$c : (X,\leq,\xi) \to (X^I,\leq,\const_{\xi})$ is given by $c(x) = \const_x$. By \ref{ccomplete}, $(C,\leq)$-${\bf poTop}$ is complete. Since $P$ has a left
adjoint (cf. \ref{adjoints}), it preserves all limits. By \ref{pullbackfibration}, the class of $(C,\leq)$-difibrations contains all isomorphisms and is closed under
base change and composition. By \ref{allobjectsfibrant}, every final morphism is a $(C,\leq)$-difibration. Using the adjunction between $P$ and the cylinder functor (cf. \ref{adjoints}) one easily sees that the homotopy lifting property is equivalent to
the defining property of $(C,\leq)$-difibrations. For a $(C,\leq)$-pospace $(X,\leq,\xi)$ the $(C,\leq)$-dimap $T : ((X^I)^I, \leq, \const_{\const_{\xi}})  \to
((X^I)^I, \leq, \const_{\const_{\xi}}) $ is given by $T(\omega)(s)(t) = \omega (t)(s)$. It remains to check P4. Let $p : (E,\leq, \varepsilon) \to
(B,\leq,\beta)$ be a $(C,\leq)$-difibration.  We have to show that the $(C,\leq)$-dimap
$$(ev_0,ev_1,p^I): (E^I,\leq,\const_{\varepsilon})\to ((E\times E)\times _{B\times B} B^I, \leq, (\varepsilon,\varepsilon, \const_{\beta}))$$ is a
$(C,\leq)$-difibration. Consider a $(C,\leq)$-dimap $f: (X,\leq,\xi) \to (E^I, \leq, \const_{\varepsilon})$ and a dimap $F : (X,\leq) \times (I,\Delta) \to
((E\times E)\times _{B\times B} B^I, \leq)$ such that $((ev_0,ev_1,p^I) \circ f)(x) = F(x,0)$ and $F(\xi(c),t) =  (\varepsilon(c),\varepsilon(c), \const_{\beta (c)})$.
Write $F = (F_0,F_1,F_2)$ and consider the following commutative diagram of spaces where $j$ is the obvious inclusion and
$\phi$ and $G$ are given by
$\phi(x,t,0) = f(x)(t)$, $\phi(x,0,s) = F_0(x,s)$, $\phi(x,1,s) = F_1(x,s)$, and $G(x,t,s) = F_2(x,s)(t)$:
$$\xymatrix{
X \times (I\times \{0\}\cup \{0,1\}\times I) \ar[rr]^-{\phi}_{} \ar[d]^{}_{id_X\times j}
& & E \ar [d]^{p}_{}
\\
X \times I\times I \ar[rr]^{}_-{G}
& & B.
}$$
Let $(x,t,s), (x',t',s') \in  X \times (I\times \{0\}\cup \{0,1\}\times I)$ such that  $(x,t,s) \leq (x',t',s')$ in  $(X,\leq) \times (I\times \{0\}\cup \{0,1\}\times I,\Delta)$.
Then $x \leq x'$, $t = t'$, and $s=s'$. It follows that $s = 0 \Rightarrow s' =0$, $t = 0 \Rightarrow t' = 0$, and $t= 1 \Rightarrow t' = 1$. Since $f$ and $F$ are
dimaps, we obtain that $\phi(x,t,s) \leq \phi(x',t',s')$ and hence that $\phi$ is a dimap $(X,\leq) \times (I\times \{0\}\cup \{0,1\}\times I,\Delta) \to (E,\leq)$.
Moreover, $\phi(\xi(c),t,s) = \varepsilon(c)$. Since $F_2$ is a dimap, $G$ is a dimap $(X,\leq) \times (I\times I,\Delta) \to (B,\leq)$. Moreover,
$G(\xi(c),t,s) = \beta(c)$. Since $j$ is a trivial cofibration in $\bf Top$, there exists, by \ref{RLPBox}, a dimap
$H : (X,\leq) \times (I\times I,\Delta) \to (E,\leq)$ such that $p \circ H = G$, $H \circ (id_X \times j) = \phi$, and $H(\xi(c), t,s) = \varepsilon (c)$. Consider the dimap $\lambda : (X,\leq) \times (I,\Delta) \to (E^I,\leq)$ defined by
$\lambda (x,s)(t) = H(x,t,s)$. We have $(ev_0,ev_1,p^I) \circ \lambda = F$, $\lambda (x,0) = f(x)$, and $\lambda (\xi(c),s) = \const_{\varepsilon(c)}$. This shows that $(ev_0,ev_1,p^I)$ is a $(C,\leq)$-difibration.
\findem

\begin{defin}\rm \cite[I.3a]{Baues} 
Let $\bf C$ be a P-category. Two morphisms $f,g: X \to Y$ are said to be \textit{homotopic}, $f \simeq g$, if  there exists a morphism $h : X \to Y^I$ such that $q_0\circ h = f$ and $q_1\circ h = g$. A morphism $f : X \to Y$ is called a \textit{homotopy equivalence} if there exists a morphism $g : Y \to X$ such that $g\circ f \simeq id_X$ and $f\circ g \simeq id_Y$.  
\end{defin}

By \ref{pathhomotopic}, two $(C,\leq)$-dimaps are homotopic in the P-category $(C,\leq)$-${\bf poTop}$ if and only if they are dihomotopic relative to $(C,\leq)$. A $(C,\leq)$-dimap is a homotopy equivalence in the P-category $(C,\leq)$-${\bf poTop}$ if and only if it is a dihomotopy equivalence relative to $(C,\leq)$.

The main result of the homotopy theory of a P-category is that it is a fibration category. There is an 
extensive homotopy theory available for fibration categories (cf. \cite{Baues}).

\begin{defin} \cite[I.1a]{Baues} \rm
A category $\bf F$ equipped with two classes of morphisms, \textit{weak equivalences} and \textit{fibrations}, is a \textit{fibration category} if it has a final object $\ast$ and if the following axioms are satisfied:
\begin{itemize}
\item[F1] An isomorphism is a \textit{trivial fibration}, i.e., a morphism which is both a fibration and a weak equivalence.
The composite of two fibrations is a fibration. If two of the morphisms $f: X \to Y$, $g:Y \to Z$, and $g\circ f: X \to Z$ are weak equivalences, so is the third.
\item[F2] The pullback of two morphisms one of which is a fibration exists. The fibrations are stable under base change. The base extension of a weak equivalence along a fibration is a weak equivalence.
\item[F3] Every morphism $f$ admits a factorization $f = p\circ j$ where $p$ is a fibration and $j$ is weak equivalence.
\item[F4] For each object $X$ there exists a trivial fibration $Y \to X$ such that $Y$ is \textit{cofibrant}, i.e., every trivial fibration $E \to Y$ admits a section.
\end{itemize}
An object $X$ is said to be \textit{$\ast$-fibrant} if the final morphism $X \to \ast$ is a fibration.
\end{defin}

Note that in \cite{Baues} a fibration category is not required to have a final object.

\begin{theor} \cite[I.3a.4]{Baues} \label{Pfib}
Let $\bf C$ be a P-category. Then $\bf C$ is a fibration category. 
The fibrations are those of the P-category $\bf C$ and the weak equivalences are the homotopy equivalences. All objects are $\ast$-fibrant and cofibrant. 
\end{theor}

As a consequence of \ref{Pfib} and \ref{PpoTop} we obtain

\begin{theor} \label{fibmain}
Let $(C,\leq)$ be a pospace. The category $(C,\leq)$-${\bf poTop}$ of $(C,\leq)$-pospaces is a fibration category. The weak equivalences are the dihomotopy equivalences relative to $(C,\leq)$
and the fibrations are the $(C,\leq)$-difibrations. All objects are $\ast$-fibrant and cofibrant.
\end{theor}

\begin{rem} \rm \label{factorization} Let $\bf C$ be a P-category. By F3, every morphism $f: X \to Y$ admits a factorization $f = p\circ j$ where $p : W \to Y$ is a fibration and $j : X \to W$ is a homotopy equivalence. By \cite[I.3]{Baues}, the object $W$ can be chosen to be the \textit{mapping path object} of $f$, i.e., the fibered product $W = X \times _YY^I$ of the morphisms $f$ and $q_0$. The homotopy equivalence $j$ is then the morphism $(id_X, c\circ f)$ and the fibration $p$ is the composite $q_1 \circ pr_{Y^I}$. 
\end{rem}

\begin{defin}\rm \cite[I.1a]{Baues} 
Let $\bf F$ be a fibration category, $p: E \to B$ be a fibration, and $X$ be a cofibrant object. Two morphisms $f, g: X \to E$ are said to be \textit{homotopic over $B$} if for some factorization of the morphism $(id_E,id_E) : E \to E\times _BE$ in a weak equivalence $E\to P$
and a fibration $q: P \to E\times_BE$ there exists a morphism $h : X \to P$ such that $q \circ h = (f,g)$. Two morphisms $f,g : X \to Y$ from a cofibrant object to a $\ast$-fibrant object are said to be \textit{homotopic}, $f \simeq g$, if they are homotopic over the final object. A morphism $f : X \to Y$ between $\ast$-fibrant and cofibrant objects is said to be a \textit{homotopy equivalence} if there exists a morphism $g : Y \to X$ such that $g\circ f \simeq id_X$ and $f\circ g \simeq id_Y$.
\end{defin}
 
\begin{prop} \label{FhoverB} 
Let $\bf C$ be a P-category and $p: E \to B$ be a fibration.  Two morphisms $f, g : X \to E$ satisfying $p\circ f = p\circ g$ are homotopic over $B$ in the fibration category $\bf C$ if and only if there exists a morphism $h : X \to E^I$ such that $q_0\circ h = f$, $q_1 \circ h = g$ and $p^I\circ h = c\circ p \circ f = c\circ p \circ g$. 
\end{prop}

\dem
We first construct a factorization of the morphism $(id_E,id_E) : E \to E\times _BE$ in a weak equivalence and a fibration. Consider the following commutative diagram:
$$\xymatrix{
B \ar[r]^{c}_{} \ar[d]^{}_{id_B} 
& B^I\ar [d]^{(q_0,q_1)}_{}
& E^I\ar[l]_{p^I}_{} \ar[d]^{(q_0,q_1)}_{} 
\\ 
B \ar[r]^{}_-{(id_B,id_B)} 
& B\times B
& E\times E. \ar[l]^{p\times p}
}$$
By the dual of the gluing lemma \cite[II.1.2]{Baues}, $p\times p$ is a fibration. By P4 and P2, $p^I$ is a fibration. We can therefore form the pullbacks of the lines of the above diagram. Applying P4 to the final morphisms $E\to \ast$ and $B\to \ast$ we obtain that the vertical morphisms are fibrations. Since, again by P4, $(q_0,q_1,p^I) : E^I \to (E\times E)\times _{B\times B}B^I$ is a fibration, we may apply the dual of the gluing lemma to deduce that the morphism $$id_B \times _{(q_0,q_1)}(q_0,q_1) : B\times _{B^I}E^I \to B\times _{B\times B}(E\times E) = E\times _BE$$
is a fibration. Consider now the following commutative diagram:
$$\xymatrix{
B \ar[r]^{id_B}_{} \ar[d]^{}_{id_B} 
& B\ar [d]^{c}_{}
& E \ar[l]^{}_{p} \ar[d]^{c}_{} 
\\ 
B \ar[r]^{}_{c} 
& B^I
& E^I. \ar[l]^{p^I}_{}
}$$
For any object $X$ the natural morphism $c : X \to X^I$ is the weak equivalence of the mapping path factorization of $id_X$. Therefore the vertical morphisms in the diagram are weak equivalences and we may apply the dual of the gluing lemma to deduce that the morphism
$$id_B\times _{c}c : E = B\times _BE \to B\times_{B^I}E^I$$
is a weak equivalence. The composite 
$$(id_B \times _{(q_0,q_1)}(q_0,q_1)) \circ (id_B\times _{c}c) : B\times _BE \to B\times _{B\times B}(E\times E)$$
is precisely the morphism $(id_E,id_E) : E \to E\times _BE$. 

Let $f,g : X \to E$ be two morphisms such that $p\circ f = p\circ g$. It follows from \cite[II.2.2]{Baues} that $f$ and $g$ are homotopic over $B$ if and only if there exists a morphism $H : X \to B \times_{B^I}E^I$ such that the following diagram is commutative: 
$$\xymatrix{
X \ar[r]^-{H}_{} \ar[d]^{}_{(f,g)} 
& B\times _{B^I}E^I \ar [d]^{id_B \times _{(q_0,q_1)}(q_0,q_1)} 
\\ 
E\times _BE \ar[r]^{}_-{id} 
& B\times _{B\times B}(E\times E).
}$$
This is the case if and only if there exists a morphism $h : X \to E^I$ such that $q_0\circ h = f$, $q_1 \circ h = g$ and $p^I\circ h = c\circ p \circ f = c\circ p \circ g$. The correspondance between $H$ and $h$ is given by $H = (p\circ f,h)$.
\findem

\begin{prop}  \label{hoverB}
Let $(C,\leq)$ be a pospace, $p : (E,\leq, \varepsilon) \to (B,\leq, \beta)$ be a $(C,\leq)$-difibration, and
$f, g : (X,\leq, \xi) \to (E,\leq, \varepsilon)$ be two $(C,\leq)$-dimaps such that \linebreak $p\circ f = p\circ \nolinebreak g$. Then
$f$ and $g$ are homotopic over $(B,\leq)$ in the fibration category $(C,\leq)$-${\bf poTop}$ if and only if there exists a
dihomotopy relative to $(C,\leq)$ $H : (X,\leq) \times (I,\Delta) \to (E,\leq)$ from $f$ to $g$ such that $p(H(x,s)) = p(f(x)) = p(g(x))$ $(x \in X$, $s\in I)$. In particular, two $(C,\leq)$-dimaps are homotopic in the fibration category $(C,\leq)$-${\bf poTop}$ if and only if they are dihomotopic relative to $(C,\leq)$ and a $(C,\leq)$-dimap is a homotopy equivalence in the fibration category $(C,\leq)$-$\bf poTop$ if and only if it is a dihomotopy equivalence relative to $(C,\leq)$. 
\end{prop}

\dem 
This follows from \ref{FhoverB}. 
\findem

\section{Cofibrations in a fibration category} \label{cofibrations}

Throughout this section we work in a fibration category $\bf F$. We suppose that all objects are cofibrant and $\ast$-fibrant and that a morphism is a weak equivalence if and only if it is a homotopy equivalence.

\begin{defin}\rm
A \textit{cofibration} is a morphism having the left lifting property with respect to the  trivial fibrations.
\end{defin}

For general reasons we have 

\begin{prop} \label{pushoutcofibration}
The class of cofibrations is closed under composition, retracts, and cobase change. Every isomorphism is a cofibration.
\end{prop}

\begin{defin}\rm
A \textit{trivial cofibration} is a cofibration which is also a weak equivalence.
\end{defin}

\begin{prop} \label{trivdicof}
A morphism is a trivial cofibration if and only if it has the left lifting property with respect to the fibrations.
\end{prop}

\dem
Let $i :A \to X$ be a morphism. Suppose first that $i$ has the left lifting property with respect to the fibrations. Then $i$ is a cofibration. Choose a factorization $i = p\circ h$ where $h : A \to E$ is a weak equivalence and $p : E \to X$ is a fibration. Thanks to our hypothesis there exists a morphism $\lambda : X\to E$ such that $p\circ \lambda = id_{X}$ and $\lambda \circ i = h$. The weak equivalence $h$ is a homotopy equivalence. Let $g$ be a homotopy inverse of $h$. We have $g\circ \lambda \circ i = g \circ h \simeq id_A$ and $i \circ g\circ \lambda = p \circ h \circ g \circ \lambda \simeq p\circ \lambda = id_X$. Thus, $i$ is a homotopy equivalence. Thanks to our general hypothesis, $i$ is a weak equivalence. 

Now suppose that $i$ is a trivial cofibration and consider a commutative diagram
$$\xymatrix{
A \ar[r]^{f}_{} \ar[d]_{i}_{} 
& E \ar [d]^{p} 
\\ 
X \ar[r]^{}_{g} 
& B
}$$
where $p$ is a fibration. Form the pullback 
$$\xymatrix{
X\times_BE \ar[r]^-{pr_E}_{} \ar[d]_-{pr_X}_{} 
& E \ar [d]^{p} 
\\ 
X \ar[r]^{}_{g} 
& B
}$$
and choose a factorization of the induced morphism $(i,f) : A \to X\times_BE$ in a weak equivalence $h : A \to Y$ and a fibration $q : Y \to X\times_BE$. Since fibrations are stable under base change and composition, $pr_X \circ q$ is a fibration. By F1, $pr_X \circ q$ is a weak equivalence. Consider the following commutative diagram:
$$\xymatrix{
A \ar[r]^-{h}_{} \ar[d]^{}_{i} 
& Y\ar [d]^-{pr_X\circ q}_{} 
\\ 
X \ar[r]^{}_{id_X} 
& X.
}$$
Since $i$ is a cofibration, there exists a morphism $\lambda : X\to Y$ such that $\lambda \circ i = h$ and $pr_X\circ q \circ \lambda = id_X$. We have $(pr_E\circ q\circ \lambda) \circ i = pr_E\circ q\circ h = f$
and $p\circ (pr_E\circ q\circ \lambda) = g \circ pr_X \circ q \circ \lambda = g$. This shows that $i$ has the required lifting property.
\findem

\begin{cor} \label{pushouttrivdicof}
The class of trivial cofibrations is closed under cobase change, composition, and retracts. Every isomorphism is a trivial cofibration. 
\end{cor}

\begin{prop} \label{allobjectsemptycof}
Suppose that $\bf F$ has an initial object $\emptyset$. For each object $X$ the initial morphism $\emptyset  \to X$ is a cofibration.
\end{prop}

\dem
Let $X$ be any object. Consider a commutative diagram 
$$\xymatrix{
\emptyset \ar[r]^{}_{} \ar[d]^{}_{} 
& E \ar [d]^{p} 
\\ 
X \ar[r]^{}_{f} 
& B
}$$
where $p$ is a trivial fibration. Form the pullback 
$$\xymatrix{
X\times_BE \ar[r]^-{pr_E}_{} \ar[d]_-{pr_X}_{} 
& E \ar [d]^{p} 
\\ 
X \ar[r]^{}_{f} 
& B.
}$$
By \cite[I.1.4]{Baues}, $pr_X$ is a trivial fibration. Since $X$ is cofibrant, $pr_X$ admits a section $s$. We have $p\circ pr_E \circ s = f$. This implies that the morphism $\emptyset  \to X$ is a cofibration.
\findem

The following concept of a cofibration category is formally dual to the one of a fibration category. For every result on fibration categories there is a dual result on cofibration categories and vice versa. 

\begin{defin}\rm \cite[I.1]{Baues}
A category $\bf C$ equipped with two classes of morphisms, \textit{weak equivalences} and \textit{cofibrations}, is a \textit{cofibration category} if it has an initial object $\emptyset$ and if the following axioms are satisfied:
\begin{itemize}
\item[C1] An isomorphism is a \textit{trivial cofibration}, i.e., a morphism which is both a cofibration and a weak equivalence.
The composite of two cofibrations is a cofibration. If two of the morphisms $f: X \to Y$, $g:Y \to Z$, and $g\circ f: X \to Z$ are weak equivalences, so is the third.
\item[C2] The pushout of two morphisms one of which is a cofibration exists. The cofibrations are stable under cobase change. The cobase extension of a weak equivalence along a cofibration is a weak equivalence. 
\item[C3] Every morphism $f$ admits a factorization $f = r\circ i$ where $i$ is a cofibration and $r$ is weak equivalence.
\item[C4] For each object $X$ there exists a trivial cofibration $X \to Y$ such that $Y$ is \textit{fibrant}, i.e., every trivial cofibration $Y \to Z$ admits a retraction.
\end{itemize}
An object $X$ is said to be \textit{$\emptyset$-cofibrant} if the initial morphism $\emptyset \to X$ is a cofibration.
\end{defin}

\begin{theor} \label{cofibrationcategory}
Suppose that $\bf F$ has an initial object $\emptyset$, that the pushout of two morphisms one of which is a cofibration exists, and that for every object $X$ the morphism $(id_X,id_X) : X\coprod X \to X$ admits a factorization in a cofibration followed by a weak equivalence. Then $\bf F$ is a cofibration category. All objects are $\emptyset$-cofibrant and fibrant.
\end{theor}

\dem
C1 follows from \ref{pushouttrivdicof}, \ref{pushoutcofibration}, and F1. By \ref{allobjectsemptycof}, all objects are $\emptyset$-cofibrant. By \ref{trivdicof}, all objects are fibrant and C4 holds. We next prove C3. Let $f : X \to Y$ be a morphism. Choose a factorization of the morphism $(id_X,id_X) : X \coprod X \to X$ in a cofibration $j : X \coprod X \to IX$ and a weak equivalence $p: IX \to X$. Note that $X \coprod X$ exists since all objects are $\emptyset$-cofibrant. Denote the canonical morphisms $X \to X \coprod X$ by $i_0$ and $i_1$. By \ref{pushoutcofibration}, $i_0$ and $i_1$ are cofibrations. By \ref{pushoutcofibration} and C1, both composites $j\circ i_0$ and $j\circ i_1$ are trivial cofibrations. Form the pushout 
$$\xymatrix{
X \ar[r]^-{f}_{} \ar[d]^{}_-{j\circ i_1} 
& Y \ar [d]^{\iota} 
\\ 
IX  \ar[r]^{}_-{\bar f} 
& Z.
}$$
By \ref{pushouttrivdicof}, $\iota$ is a trivial cofibration. Let $r : Z \to Y$ be the morphism induced by the morphisms $f\circ p: IX  \to Y$ and $id_{Y}$. Since $r \circ \iota = id_{Y}$ and $\iota$ and $id_{Y}$ are weak equivalences, $r$ is a weak equivalence. Let $i$ be the composite of the morphisms $j\circ i_0 : X \to IX$ and $\bar f : IX \to Z$. Consider the following pushout diagram:
$$\xymatrix{
X \coprod X \ar[r]^-{id_X \coprod f}_{} \ar[d]^{}_-{j} 
& X \coprod Y \ar [d]^{(i,\iota)} 
\\ 
IX \ar[r]^{}_-{\bar f} 
& Z.
}$$
Since $j$ is a cofibration, $(i,\iota)$ is a cofibration. Since the initial morphism $\emptyset \to Y$ is a cofibration and cofibrations are closed under cobase change (cf. \ref{pushoutcofibration}), the canonical morphism $\phi : X \to X \coprod Y$ is a cofibration. Since the composite of cofibrations is a cofibration (cf. \ref{pushoutcofibration}), $i = (i, \iota) \circ \phi$ is a cofibration. We have $r\circ i = r\circ \bar f \circ j\circ i_0 = f\circ p \circ j\circ i_0 =f$. This shows that C3 holds. C2 follows from \ref{pushoutcofibration}, \ref{pushouttrivdicof}, C1, C3, and \cite[I.1.4]{Baues}.
\findem

\begin{rem}\rm
The factorization $f = r \circ i$ constructed in the above proof is the \textit{mapping cylinder factorization} of $f$ which is dual to the mapping path factorization of \ref{factorization}.
\end{rem}

\begin{defin}\rm \cite[I.1]{Baues}
Let $\bf C$ be a cofibration category.  
Two morphisms $f,g : X \to Y$ from a $\emptyset$-cofibrant object to a fibrant object are said to be \textit{homotopic} if for some factorization of the morphism $(id_X,id_X) : X \coprod X \to X$ in a cofibration $i : X\coprod X \to IX$ and a weak equivalence $r: IX \to X$ there exists a morphism $H: IX \to Y$ such that $H\circ i = (f,g): X\coprod X \to Y$.
\end{defin}

\begin{prop} \label{ho=ho}
Two morphisms of $\bf F$ are homotopic in the cofibration category $\bf F$ if and only if they are homotopic in the fibration category $\bf F$. 
\end{prop}

\dem
Let $\simeq$ denote the homotopy relation of the fibration category $\bf F$ and $\sim$ denote the homotopy relation of the cofibration category $\bf F$. Both relations are natural equivalence relations (c.f. \cite[II.3.2]{Baues}) and we can form the quotient categories ${\bf F}/\simeq$ and ${\bf F}/\sim$. By \cite[II.3.6]{Baues}, both quotient categories have the universal property of the localization of $\bf F$ with respect to the weak equivalences. 
This implies that there is an isomorphism of categories ${\bf F}/\simeq \to {\bf F}/\sim$ which is the identity on objects and which sends the $\simeq$-class of a morphism to its $\sim$-class. 
The result follows.
\findem

\section{$(C,\leq)$-Dicofibrations} \label{dicofibrations}

\begin{defin}\rm
Let $(C,\leq)$ be a pospace. A \textit{$(C,\leq)$-dicofibration} is a $(C,\leq)$-dimap having the left lifting property with respect to the trivial $(C,\leq)$-difibrations. If $C = \emptyset$ we simply talk of dicofibrations.
\end{defin}

\begin{theor} \label{pocofcat}
Let $(C,\leq)$ be a pospace. The category $(C,\leq)$-${\bf poTop}$ is a cofibration category. The cofibrations are the $(C,\leq)$-dicofibrations and the weak equivalences are the dihomotopy equivalences relative to $(C,\leq)$. All objects are fibrant and $\emptyset$-cofibrant. Two $(C,\leq)$-dimaps are homotopic in the cofibration category $(C,\leq)$-${\bf poTop}$ if and only if they are dihomotopic relative to $(C,\leq)$.   
\end{theor}

\dem
Thanks to \ref{ccomplete}, \ref{fibmain}, \ref{hoverB}, \ref{cofibrationcategory}, and \ref{ho=ho} it is enough to show that for every $(C,\leq\nolinebreak)$-pospace $(X,\leq,\xi)$ the $(C,\leq)$-dimap
$(id_X,id_X) : (X,\leq,\xi) \coprod (X,\leq,\xi) \to (X,\leq,\xi)$ admits a factorization in a $(C,\leq)$-dicofibration and a dihomotopy
equivalence relative to $(C,\leq)$. Let $(X,\leq,\xi)$ be a $(C,\leq)$-pospace. We have
$$(X,\leq,\xi) \coprod (X,\leq,\xi) = (X,\leq,\xi) \Box_{(C,\leq)} \{0,1\}$$
and $(id_X,id_X)$ is the $(C,\leq)$-dimap $(X,\leq,\xi) \Box_{(C,\leq)} \{0,1\} \to (X,\leq,\xi)$, $[x,t] \mapsto x$. Let $\iota : \{0,1\} \hookrightarrow I$
be the inclusion. We show that $$(X,\leq,\xi) \Box_{(C,\leq)} \iota : (X,\leq,\xi) \Box_{(C,\leq)} \{0,1\} \to (X,\leq,\xi) \Box_{(C,\leq)} I$$ is a $(C,\leq)$-dicofibration and
that the projection $r : (X,\leq,\xi) \Box_{(C,\leq)} I \to (X,\leq,\xi)$, $r([x,t]) = x$ is a dihomotopy equivalence relative to $(C,\leq)$.
The $(C,\leq)$-dimap $\sigma : (X,\leq,\xi) \to (X,\leq,\xi) \Box_{(C,\leq)} I$ given by $\sigma (x) = [x,0]$ is a dihomotopy inverse relative to
$(C,\leq)$ of $r$. Indeed, $r \circ \sigma = id_X$ and a dihomotopy relative to $(C,\leq)$ from $\sigma \circ r$  to $id_{X\Box_C I}$ is given
by $F([x,t],s) = [x,st]$.

We now show that  $(X,\leq,\xi) \Box_{(C,\leq)} \iota$ is a $(C,\leq)$-dicofibration. Consider a commutative diagram of $(C,\leq)$-pospaces
$$\xymatrix{
(X,\leq,\xi) \Box_{(C,\leq)} \{0,1\} \ar[rr]^-{f}_{} \ar[d]^{}_{X\Box_C \iota}
& & (E,\leq, \varepsilon) \ar [d]^{p}_{}
\\
(X,\leq,\xi) \Box_{(C,\leq)} I \ar[rr]^{}_-{g}
& & (B,\leq,\beta)
}$$
where $p$ is a trivial $(C,\leq)$-difibration. By the dual of the lifting lemma \cite[II.1.11]{Baues}, there exists a section $s$ of $p$ such that $s\circ p$ is
homotopic to $id_{(E,\leq,\varepsilon)}$ over \linebreak $(B,\leq,\beta)$ in the fibration category $(C,\leq)$-${\bf poTop}$. By \ref{hoverB},
there exists a dihomotopy relative to $(C,\leq)$ $H : (E,\leq) \times (I,\Delta) \to (E,\leq)$ from $s\circ p$ to $id_{(E,\leq,\varepsilon)}$
such that $p(H(x,\tau)) = p(x)$. Consider the following commutative diagram of spaces  where $j$ is the obvious inclusion and
$\phi$ and $\Phi$ are
given by $\phi(x,t,0) = s(g([x,t]))$, $\phi(x,0,\tau) = H(f([x,0]),\tau)$, $\phi(x,1,\tau ) = H(f([x,1]),\tau)$, and $\Phi(x,t,\tau ) = g([x,t])$:
$$\xymatrix{
X \times (I\times \{0\}\cup \{0,1\}\times I) \ar[rr]^-{\phi}_{} \ar[d]^{}_{id_X\times j}
& & E \ar [d]^{p}_{}
\\
X \times I\times I \ar[rr]^{}_-{\Phi}
& & B.
}$$
Since $g$ is a dimap, $\Phi$ is a dimap $(X,\leq) \times (I\times I,\Delta) \to (B,\leq)$. Moreover,
$\Phi (\xi(c),t,\tau) = g(\bar \xi (c)) = \beta (c)$. Let $(x,t,\tau)$, $(x',t',\tau') \in X \times (I\times \{0\}\cup \{0,1\}\times I)$ such that
$(x,t,\tau) \leq (x',t',\tau')$ in $(X,\leq) \times (I\times \{0\}\cup \{0,1\}\times I,\Delta)$. Then $x \leq x'$, $t = t'$, and $\tau = \tau'$. It follows that $t = 0 \Rightarrow t' = 0$, $t = 1 \Rightarrow t' = 1$, and $\tau = 0 \Rightarrow \tau ' = 0$. We obtain that $\phi(x,t,\tau) \leq \phi(x',t',\tau')$. Thus $\phi$ is a dimap \linebreak $(X,\leq) \times (I\times \{0\}\cup \{0,1\}\times I,\Delta) \to (E,\leq)$.
Moreover, $\phi (\xi(c),t,\tau) = \varepsilon (c)$. \linebreak Since $j$  is a trivial cofibration of spaces,
there exists, by \ref{RLPBox}, a dimap \linebreak 
$G:(X,\leq) \times (I \times I,\Delta) \to (E,\leq)$ such that $G \circ (id_X\times j) = \phi$, $p\circ G = \Phi$, and $G(\xi(c),t,\tau) = \varepsilon (c)$.
Let
$\lambda : (X,\leq,\xi) \Box_{(C,\leq)} I \to (E, \leq,\varepsilon)$ be the $(C,\leq)$-dimap given by
$\lambda ([x, t]) = G(x,t,1)$. We have
$\lambda ([x,0]) = G(x,0,1) = \phi(x,0,1) = H(f([x,0]),1) = f([x,0]),$
$\lambda ([x,1]) = G(x,1,1) = \phi(x,1,1)  = H(f([x,1]),1) = f([x,1]),$
and $(p\circ \lambda) ([x,t]) = (p \circ G) (x,t,1) = \Phi(x,t,1) =g([x,t])$. This shows that $(X,\leq,\xi)\Box_{(C,\leq)} \iota$ is a $(C,\leq)$-dicofibration.
\findem

\section{The closed model category of pospaces} \label{modelcategory}
In this section we show that absolute pospaces form a closed model category. Closed model categories have been introduced by D. Quillen \cite{Quillen}. The definition below can for instance be found in the book \cite{GoerssJardine} by P. Goerss and J.F. Jardine.

\begin{defin} \rm
A category $\bf C$ equipped with three classes of morphisms, \textit{weak equivalences}, \textit{fibrations} and \textit{cofibrations}, is a \textit{closed model category} if the following axioms are satisfied:
\begin{itemize}
\item[CM1] The category $\bf C$ is finitely complete and cocomplete.
\item[CM2] If two of the morphisms $f: X \to Y$, $g:Y \to Z$, and $g\circ f: X \to Z$ are weak equivalences, so is the third.
\item[CM3] Weak equivalences, fibrations, and cofibrations are closed under retracts.
\item[CM4] The fibrations have the right lifting property with respect to the \textit{trivial cofibrations}, i.e. the morphisms which are both cofibrations and weak equivalences. The cofibrations have the left lifting property with respect to the \textit{trivial fibrations}, i.e. the morphisms which are both fibrations and weak equivalences. 
\item[CM5] Every morphism $f$ admits a factorization $f = p\circ i$ where $p$ is a fibration and $i$ is a trivial cofibration. Every morphism $f$ admits a factorization $f = q\circ j$ where $j$ is a cofibration and $q$ is a trivial fibration. 
\end{itemize}
An object is said to be \textit{fibrant} if its final morphism is a fibration. An object is said to be \textit{cofibrant} if its initial morphism is a cofibration. Two morphisms $f,g : X \to Y$ from a cofibrant object to a fibrant object are said to be \textit{homotopic} if for some factorization of the morphism $(id_X,id_X) : X \coprod X \to X$ in a cofibration $i : X\coprod X \to IX$ and a weak equivalence $r: IX \to X$ there exists a morphism $H: IX \to Y$ such that $H\circ i = (f,g): X\coprod X \to Y$.
\end{defin}

\begin{lem} \label{RLP}
Let $i : (A,\leq) \to (X,\leq)$ be an inclusion of pospaces such that there exist a dihomotopy $H: (X,\leq) \times (I,\Delta ) \to (X,\leq)$ and a dimap $\phi : (X,\leq) \to (I,\Delta)$ such that $A = \phi ^{-1}(0)$, $H(x,1) = x$ $(x\in X)$, $H(a,t) = a$ $(a\in A$, $t\in I)$, and $H(x,0) \in A$ $(x\in X)$. Then $i$ is a trivial dicofibration.
\end{lem}

\dem
The proof is similar to the one of \ref{RLPBox}. Consider a difibration $p: (E,\leq) \to (B,\leq)$, a dimap $f : (A,\leq) \to (E,\leq\nolinebreak)$, and a dimap
$g : (X,\leq) \to (B,\leq)$ such that $g(a)= p(f(a))$. Write $r(x) = H(x,0)$. 
Consider the map $G : X \times I \to X$ defined by
$$G(x,t)= \left\{ \begin{array}{lll}
H(x, \frac{t}{\phi(x)}) & t < \phi(x), \\
x & t \geq \phi(x).
\end{array}\right. $$
As in \cite[I.7.15]{Whitehead} one shows that $G$ is continuous. We have $G(x,0) = (i\circ r)(x)$ for all $x \in X$. Let $(x,t), (x',t') \in X \times I$ such that $(x,t) \leq (x',t')$. Then $x\leq x'$ and $t = t'$. Since $\phi$ is a dimap, $\phi(x) \leq \phi(x')$ and hence $\phi(x) = \phi(x')$. It follows that $G(x,t) \leq G(x',t')$ and hence that $G$ is a dimap. Consider the following commutative diagram of
pospaces:
$$\xymatrix{
(X,\leq) \ar[rr]^-{f\circ  r}_{} \ar[d]^{}_{i_0}
& &(E,\leq)\ar [d]^{p}_{}
\\
(X,\leq)\times (I,\Delta) \ar[rr]^{}_-{g\circ G}
& & (B,\leq).
}$$
Since $p$ is a difibration, there exists a dimap $F : (X,\leq) \times (I,\Delta) \to (E,\leq)$ such that
$F \circ i_0 = f\circ r$ and  $p\circ F = g\circ G$. Consider the dimap
$\lambda : (X,\leq)  \to (E,\leq)$ defined by $\lambda (x) = F(x,\phi(x))$. We have
$(p\circ \lambda) (x) = p(F(x,\phi(x))) = g(G(x,\phi(x))) = g(x)$
and
$\lambda (a) = F(a,\phi(a)) = F(a,0) = f(r(a)) = f(a).$
By \ref{trivdicof}, this shows that $i$ is a trivial dicofibration. 
\findem

\begin{theor} \label{model}
The category $\bf poTop$ of pospaces is a closed model category where weak equivalences are dihomotopy equivalences, fibrations are difibrations, and cofibrations are dicofibrations. All pospaces are fibrant and cofibrant. Two dimaps are homotopic in the closed model category ${\bf poTop}$ if and only if they are dihomotopic.  
\end{theor}

\dem
By \ref{complete}, $\bf poTop$ is complete and cocomplete. CM2 is part of \ref{2=3}. CM3 follows from \ref{2=3}, \ref{pullbackfibration}, and \ref{pushoutcofibration}. CM4 follows from the definition of dicofibrations and \ref{trivdicof}. We now show CM5. Let $f : (X,\leq) \to (Y,\leq)$ be a dimap. We show first that $f$ admits a factorization $f = p\circ i$ where $p$ is a fibration and $i$ is a trivial cofibration. We proceed as in \cite{Strom}. Consider the mapping path factorization $f = q \circ j$ of \ref{factorization}; $j : (X,\leq) \to (X\times_YY^I,\leq\nolinebreak )$ is the dihomotopy equivalence given by $j(x) = (x,\const_{f(x)})$ and $q$ is the difibration $ (X\times_YY^I,\leq) \to (Y,\leq)$, $(x,\omega ) \mapsto \omega (1)$. Let $(E,\leq)$ be the subpospace of $(X\times_YY^I,\leq) \times (I,\Delta)$ defined by $E = j(X) \times I \cup (X\times_YY^I) \times ]0,1]$. Let \linebreak $i : (X,\leq) \to (E,\leq)$ be the dimap defined by $i(x) = (j(x),0)$. Let $p : (E,\leq) \to (Y,\leq)$ be the composite
$$(E,\leq) \hookrightarrow (X\times_YY^I,\leq) \times (I,\Delta) \stackrel{pr_{X\times_YY^I}}{\longrightarrow} (X\times_YY^I,\leq) \stackrel{q}{\to} (Y,\leq).$$
We have $p \circ i = f$. We show that $i$ is a trivial dicofibration. Consider the dimap $H: (E,\leq)\times (I,\Delta) \to (E,\leq)$ defined by $H(x,\omega,s,t) = (x,\omega _{t},st)$ where $\omega _{t}(\tau) = \omega (t\tau)$. We have $H(x,\omega,s,1) = (x,\omega,s)$, $H(x,\const_{f(x)},0,t) = (x,\const_{f(x)},0)$, and $H(x,\omega,s,0) = (x,\const_{f(x)},0) \in i(X)$. Moreover, $i(X) = \phi ^{-1}(0)$ where $\phi : (E,\leq\nolinebreak) \to (I,\Delta)$ is the dimap $(x,\omega, t) \mapsto t$. By \ref{RLP}, the inclusion $(i(X),\leq) \hookrightarrow (E,\leq)$ is a trivial dicofibration. The dimap $i:(X,\leq) \to (i(X),\leq)$ is an isomorphism of pospaces; the inverse the composite $$(i(X),\leq) \hookrightarrow (E,\leq) \hookrightarrow (X\times_YY^I,\leq) \times (I,\Delta) \stackrel{pr_{X\times_YY^I}}{\longrightarrow} (X\times_YY^I,\leq) \stackrel{pr_X}{\to} (X,\leq).$$
By \ref{pushouttrivdicof}, it follows that $i$ is a trivial dicofibration. 

We now show that $p$ is a difibration. Since the composite of two difibrations is a difibration and $q$ is a difibration, it suffices to show that the composite 
$$\pi : (E,\leq) \hookrightarrow (X\times_YY^I,\leq) \times (I,\Delta) \stackrel{pr_{X\times_YY^I}}{\longrightarrow} (X\times_YY^I,\leq)$$
is a difibration. Consider a commutative diagram of pospaces
$$\xymatrix{
(Z,\leq) \ar[r]^{g}_{} \ar[d]^{}_{i_0} 
& (E,\leq) \ar [d]^{\pi}_{} 
\\ 
(Z,\leq) \times (I,\Delta) \ar[r]^{}_{G} 
& (X\times_YY^I,\leq).
}$$
Consider the continuous map $F : Z \times I \to E$ defined by $$F(z,t) = (G(z,t),t + (1-t)\phi(g(z))).$$
If $t + (1-t)\phi(g(z)) = 0$ then $t =0$ and $\phi (g(z))=0$. Therefore $g(z) \in i(X)$ and $$F(z,t) = (G(z,0),0) = (\pi (g(z)),0) = g(z) \in i(X) \subset E.$$
This shows that $F$ is well-defined. Let $(z,t) \leq (z',t') \in Z\times I$. Then $z \leq z'$ and $t = t'$. Since $\phi$ and $g$ are dimaps, $\phi(g(z)) \leq \phi(g(z'))$, i.e., $\phi(g(z)) = \phi(g(z'))$. It follows that $F(z,t) \leq F(z',t')$ and hence that $F$ is a dimap. We have
$$F(z,0) = (G(z,0),\phi(g(z))) = (\pi(g(z)),\phi(g(z))) = g(z)$$
and $\pi \circ F = G$. It follows that $\pi $ is a difibration. This terminates the proof of the first part of CM5. 
For the second part of CM5 we use \ref{pocofcat} to obtain a factorization $f = r\circ j$ where $j$ is a dicofibration and $r$ is a dihomotopy equivalence. As we have seen, $r$ admits a factorization $r = p \circ \iota$ where $\iota$ is a trivial dicofibration and $p$ is a difibration. By CM2, $p$ is a dihomotopy equivalence. Since the composite of two dicofibrations is a dicofibration, $i = \iota \circ j$ is a dicofibration. Thus, $f = p\circ i$ is a factorization as required. This terminates the proof of CM5. The remaining statements follow from \ref{fibmain} and \ref{pocofcat}.
\findem

\vspace{1cm}
\noindent \textit{Address of the author:}
Universidade do Minho,
Centro de Matem\'atica,\linebreak
Campus de Gualtar,
4710-057 Braga,
Portugal.
{\it{E-mail:}} kahl@math.uminho.pt
%, lucile@math.uminho.pt\\

\end{document}